\documentclass[11pt,a4paper]{amsart}
\usepackage{amscd}
\theoremstyle{plain} 
\newtheorem{thm}{Theorem}[section]
\newtheorem{cor}[thm]{Corollary}
\newtheorem{lem}[thm]{Lemma}
\newtheorem{prop}[thm]{Proposition}

\theoremstyle{definition}

\newtheorem{rem}[thm]{Remark}

\theoremstyle{remark}

\numberwithin{equation}{section}

\newsavebox{\SmallMathBox}

\def\dd{\partial}

\def\Di{D\kern -.65em /}
\def\Dii{D\kern -.45em /}
\def\di{{\dd}\kern -.55em /}
\def\dii{{\dd}\kern -.40em /}

\def\zd{{\det}_{\zeta}}

\def\NN{{\bf N}}

\def\Dd{{\mathcal D}}
\def\Ee{{\mathcal E}}

\def\Kk{{\mathcal K}}

\def\Nn{{\mathcal N}}

\def\Rr{{\mathcal R}}

\def\={\cong}
\def\>{\supset}
\def\<{\subset}

\def\12{\frac{1}{2}}
\def\2{\Dd}
\def\3{\Nn}
\def\4{\Rr}
\def\6{\cup}
\def\8{\otimes}
\def\0{^{\circ}}

\def\e{\varepsilon}

\def\g{\gamma}
\def\G{\Gamma}

\def\la{\lambda}

\def\N{\NN}

\def\Si{\Sigma}

\def\z{\zeta}

\def\bC{\mathbb{C}}

\def\Id{{\rm Id}}

\def\ran{\mbox{\rm range\,}}

\def\Si{S\kern -.65em /}

\def\tensor{\otimes}

\def\Tr{\mbox{\rm Tr\,}}

\def\Del{\Delta}
\def\sql{{\lambda^2}}

\begin{document}

\title[adiabatic decomposition of the $\z$-determinant]
{Adiabatic Decomposition of the $\z$-determinant and Dirichlet to
Neumann operator}

\author{Jinsung Park}
\address{Mathematisches Institut, Universit\"at Bonn \\
Beringstrasse 1, D-53115, Bonn, Germany}
\email{jpark@math.uni-bonn.de}

\author{Krzysztof P. Wojciechowski}
\address{Department of Mathematics\\IUPUI
(Indiana/Purdue)\\
Indianapolis IN 46202--3216, U.S.A.}
\email{kwojciechowski@math.iupui.edu}

\vskip 5mm

\begin{abstract}
We discuss the adiabatic decomposition formula of the
$\z$-determ-inant of a Laplace type operator on a  closed
manifold. We also analyze the adiabatic behavior of the
$\z$-determinant of a Dirichlet to Neumann operator. This analysis
makes it possible to compare the adiabatic decomposition formula
with the Meyer-Vietoris type formula for the $\z$-determinant
proved by Burghelea, Friedlander and Kappeler. As a byproduct of
this comparison, we obtain the exact value of the local constant
which appears in their formula for the case of Dirichlet boundary
condition.
\end{abstract}

\thanks{key words: $\z$-determinant, adiabatic limit, Dirichlet to
Neumann operator. } 

\maketitle

\section{Introduction and Statement of the Results}

\bigskip

In this paper we continue our study of the adiabatic decomposition
of the $\zeta$-determinant of the Laplace type operator. In
\cite{JPKW4}, \cite{JPKW5} the decomposition formula of the
$\z$-determinant of Dirac Laplacian was given in terms of the
non-local Atiyah-Patodi-Singer boundary condition. Here we discuss
a formula which involves the \emph{Laplace type operator} and the
\emph{Dirichlet boundary condition}.

\bigskip

Let $\Delta : C^{\infty}(M,E) \to C^{\infty}(M,E)$ denote a
Laplace type operator acting on sections of a vector bundle $E$
over a closed manifold $M$ of dimension $n$. The operator $\Delta$
is a self-adjoint operator with discrete spectrum $\{\la_k\}_{k
\in \N}$. Let us decompose $M$ into two sub-manifolds $M_1$, $M_2$
with common boundary $Y$,
\begin{equation}\label{e:dec1}
M = M_1 \cup M_2 \ \ , \ \ M_1 \cap M_2 = Y = {\partial}M_1 =
{\partial}M_2 \ \, .
\end{equation}
The $\z$-function $\z_{\Del}(s)$ is defined by
$$\z_{\Del}(s) =
\sum_{\la_k\neq 0}\la_k^{-s} \ \ ,$$ which is a holomorphic
function in the half-plane $\Re(s) > \frac n2$ and extends to a
meromorphic function on the whole complex plane with $s=0$ as a
regular point. The $\z$-determinant of $\Del$ is defined by
\begin{equation}\label{e:zz1}
\log \zd\Del = - \frac d{ds} \z_{\Del}(s)\Bigr|_{s=0} \ \ .
\end{equation}
The derivative of $\z_{\Del}(s)$ at $s=0$ can be represented in
the following way
\begin{equation}\label{e:zz2}
\frac d{ds} \z_{\Del}(s)\Bigr|_{s=0} = \lim_{s \to 0}(\kappa(s) -
\frac {a'_{n/2}}s) + \g a'_{n/2} \ \, .
\end{equation}
Here $\g$ denotes Euler's constant and $a'_{n/2}:=a_{n/2}-\dim
\ker(\Del)$ where $a_{n/2}$ is constant term in the following
asymptotic expansion near $t=0$,
$$\Tr (e^{-t\Del})\ \sim \ t^{- \frac n2}
\sum_{k = 0}a_kt^{\frac k 2} \ \ .$$ The function $\kappa(s)$ is
defined as the integral
\begin{equation}\label{e:kappa}
\kappa(s) = \int_0^{\infty}t^{s-1} \big(\, \Tr(e^{-t\Del}) - \dim
\ker(\Del)\, \big)\ dt \ \,
\end{equation}
for $\Re(s) > \frac{n}{2}$ . It has a meromorphic extension to the
whole complex plane and it can be represented as
$$\kappa(s) = \frac {a'_{n/2}}s + h(s)$$
in a neighborhood of $s=0$ , where $h(s)$ is a holomorphic
function of $s$ . The value of the function $h(s)$ at $s=0$ is not
a local invariant, and this fact implies the non-locality of the
$\z$-determinant. This is the main reason, that there is no
straightforward decomposition formula for the $\z$-determinant of
the operator $\Del$ onto contributions coming from $M_1$ and $M_2$
(see \cite{JPKW2}, \cite{JPKW3} and \cite{JPKW5} for more detailed
discussion).

\bigskip
We assume that there is a bicollar neighborhood $N \cong [-1,1]
\times Y$ of $Y$ in $M$ such that the Riemannian structure on $M$
and the Hermitian structure on $E$ are products of the
corresponding structures over $[-1,1]$ and $Y$ when restricted to
$N$ . We also assume that the operator $\Del$ restricted to $N$
has the following form
\begin{equation}\label{e:pr}
\Del = -\partial_u^2 + \Del_Y \ \, .
\end{equation}
Here $u$ denotes the normal variable and $\Del_Y$ is a
$u$-independent Laplace type operator on $Y$ .

\bigskip

We replace the bicollar $N$ by $N_R = [-R,R] \times Y$ to obtain a
new closed manifold $M_R$ and extend the vector bundle $E$ to
$M_R$ in an obvious way. We use formula (\ref{e:pr}) to extend
$\Del$ to the Laplace operator $\Del_R$ on $M_R$ . We decompose
$M_{R}$ into $M_{1,R}$ and $M_{2,R}$ by cutting $M_R$ at
$\{0\}\times Y$.  We denote by $\Delta_{i,R}$ the operator
$\Delta_R|_{M_{i,R}}$ subject to the Dirichlet boundary condition.
The operator $\Delta_{i,R}$ is a self-adjoint operator with
discrete spectrum and smooth eigensections. The $\z$-determinant
of $\Delta_{i,R}$ is defined as $\zd\Del_R$ and it enjoys all the
nice properties of the $\z$-determinant of the Laplacian on a
closed manifold. The concern of this paper is to investigate the
\emph{adiabatic decomposition} of $\zd \Del_R$, that is, the limit
of
\begin{equation}\label{e:111}
\frac {\zd\Delta_R}{\zd\Delta_{1,R}{\cdot} \zd\Delta_{2,R}} \qquad
\text{as} \quad R\to\infty\ \ .
\end{equation}

\bigskip

The case of the invertible tangential operator $\Delta_Y$ was
described in  \cite{JPKW1}, \cite{JPKW2} and \cite{JPKW3}. The
invertibility assumption on $\Delta_Y$ implies that we have only
finitely many eigenvalues of $\Del_R$ converging to $0$ as
$R\to\infty$. This allows us to discard the large time
contribution to the $\z$-determinant of $\Del_R$ under the
adiabatic process and the adiabatic decomposition of the
$\z$-determinant easily follows from a standard application of the
\emph{Duhamel principle}.
\bigskip

The non-invertible case was studied in \cite{JPKW5}. The
decomposition formula introduced in \cite{JPKW5} uses
Atiyah-Patodi-Singer boundary conditions. The new feature of the
non-invertible tangential operator is the presence of infinitely
many eigenvalues approaching $0$ as $R \to \infty$. The behavior
of these eigenvalues can be understood in terms of suitable
scattering operators described in  \cite{Mu94}. We used this
description of \emph{small} eigenvalues in the proof of our
decomposition formula (see \cite{JPKW5}, see also announcement
\cite{JPKW4}). Since the presented results in \cite{JPKW4},
\cite{JPKW5} hold only for the Dirac type operator, we need some
modifications to deal with the Laplace case in this paper.

\bigskip

To avoid delicate analytical issues we make one more assumption.
Let us recall the classification of the eigenvalues of a Dirac
type operator $\Dd_R$ over $M_R$. The operator $\Dd_R$ has
finitely many eigenvalues $\{\la_k(R)\}$ , which decay
exponentially as $R\to\infty$ , meaning that there exists positive
constants $c_1$ and $c_2$ such that
$$|\la_k(R)| < c_1e^{-c_2R} \ \, .$$
We called them \emph{e-values} in \cite{JPKW5}. There are also
infinite families of eigenvalues, which decay like $R^{-1}$, of
$\Dd_R$ and the restrictions of $\Dd_R$ to $M_{i,R}$ with
generalized APS spectral boundary conditions. We called those
eigenvalues \emph{s-values} in \cite{JPKW5}. Finally, we have
infinitely many eigenvalues bounded away from $0$ . By our
definition, the set of zero eigenvalues is a subset of the set of
\emph{e-values} and it is known that the set of \emph{e-values} is
stable under the adiabatic process although the set of zero
eigenvalues is not. Up to now, no analysis has been known to deal
with \emph{e-values}. In order to avoid analytical difficulties
related to exponentially small eigenvalues, throughout this paper
we assume the following condition:
\begin{equation}\label{Condition A}
\text{There are no eigenvalues of $\Delta_R$ exponentially
decaying to $0$ as $R\to\infty$.}
\end{equation}
Hence, this condition means that all the eigenvalues of $\Del_R,
\Del_{i,R}$ converging to $0$ are \emph{s-values} decaying like
$R^{-2}$. There are many natural Laplace type operators satisfying
the condition \eqref{Condition A}. For example, let
$\Delta^k_{\rho,R}$ denote the Hodge Laplacian over $M_R$ acting
on the space of $k$-forms twisted by the flat vector bundle
defined by a unitary representation $\rho$ of $\pi_1(M_R)$. Then,
as in Section 4 of \cite{Ha98}, one can show that there are no
eigenvalues of $\Delta^k_{\rho,R}$ exponentially decaying to $0$
as $R\to\infty$ if $\Delta^k_{\rho,0}$ has no zero eigenvalues.

\bigskip

Let $M_{i,\infty}$ denote the manifold $M_i$ with the half
infinite cylinder attached and $\Del_{i,\infty}$ denote the
Laplace operators on $M_{i,\infty}$ determined by $\Del_i$ . The
operator $\Del_{i,\infty}$ defines a scattering matrix $C_i(0) :
\ker( \Del_Y) \to \ker(\Del_Y)$ , which is an involution over
$\ker(\Delta_Y)$. The following theorem is the first main result
of this paper,

\bigskip

\begin{thm}\label{t:main thm}
Let us assume that $\Del_R$ satisfies \eqref{Condition A}. Then we
have
\begin{equation}\label{e:main}
\lim_{R\to\infty} R^{h_Y} \frac{\zd \Del_R}{\zd \Del_{1,R}\cdot
\zd \Del_{2,R}}= 2^{-h_Y} \sqrt{\zd^* \Del_Y} \cdot
\det\big(\frac{\Id-C_{12}}{2}\big)\ \ \, ,
\end{equation}
where $h_Y:= \dim\ker(\Del_Y)$, $C_{12}:= C_1(0)\circ C_2(0)$ is a
unitary operator and $\zd^* \Del_Y$ denotes the $\z$-determinant
of the operator $\Del_Y$ restricted to the orthogonal complement
of $\ker(\Del_Y)$ .
\end{thm}

\bigskip

\begin{rem} The condition
\eqref{Condition A} implies that the operator $C_{12}$ is a
unitary operator with no unity eigenvalues (see Remark
\ref{r:c12}). It follows that $\det\big(\frac{\Id-C_{12}}{2}\big)$
is a positive real number. The operators $\Del_{i,R}$ are
Laplacians subject to the Dirichlet conditions so that all their
eigenvalues satisfy \eqref{Condition A} by a standard application
of the mini-max principle. The formula \eqref{e:main} in Theorem
\ref{t:main thm} has been used in \cite{BP04} where the adiabatic
surgery formula of the determinant line bundle is investigated.
The related decomposition formula for the analytic torsion was
also worked out by Hassell in \cite{Ha98}. He proved the analytic
surgery formula of the analytic torsion using the
\emph{b-calculus}. We also refer to the work of
Hassell-Mazzeo-Melrose \cite{HMaM95} where the analytic surgery
problem is investigated extensively.
\end{rem}

\bigskip

Our proof of Theorem \ref{t:main thm} is modelled on a proof given
in \cite{JPKW5}, with necessary modifications since we are dealing
with a different type of boundary conditions. The main
modification is a revised relation between \emph{s-values} and the
scattering matrix $C_i(0)$ . This is the main achievement of the
first part of this paper, which consists of the following two
sections.

\bigskip

In the second part, we study the adiabatic limit of the
$\z$-determinant of certain operator $\Rr_R$ appearing in the
formula of Burghelea, Friedlander and Kappeler \cite{BFK} ( in
short, BFK from now on ). The BFK formula can be formulated in our
situation as follows,
\begin{equation}\label{e:BFK0}
\frac{\zd \Del_R}{\zd \Del_{1,R}\cdot \zd \Del_{2,R}}\ = \ C(Y)\,
\zd \Rr_R \qquad \text{for any} \ \, R \ \ ,
\end{equation}
where $C(Y)$ is a locally computable constant and $\Rr_R$ is
defined as the sum of the Dirichlet to Neumann operators over the
decomposed manifolds $M_{i,R}$ . It is well known that $\Rr_R$ is
a nonnegative pseudo-differential operator of order $1$. In
particular, under the condition \eqref{Condition A}, $\Rr_R$ is a
positive operator for any $R$.

\bigskip

\begin{rem}
 The BFK constant $C(Y)$ is locally computable from symbols
of $\Delta_R^{-1}$ over $Y$, so that $C(Y)$ may depend on the
intrinsic data over $Y$ as well as the extrinsic data out of $Y$
like the normal derivatives of the symbol of $\Delta_R^{-1}$ at
$Y$. However, under the assumption of the product structure near
$Y$ , the constant $C(Y)$ depends on only the intrinsic data over
$Y$ , in particular $C(Y)$ does not change under the adiabatic
process.
\end{rem}

\bigskip

In Section 4, we study the adiabatic limit of $\zd\Rr_R$. Here we
consider the case of the non-invertible tangential operator
$\Del_Y$ , as a result, the adiabatic limit of $\zd \Rr_R$
contains the contribution determined by $\Delta_Y$ as well as the
scattering data. The following theorem is the main result for
this,

\bigskip

\begin{thm}\label{t:DN0}
Let us assume \eqref{Condition A}. Then we have the following
formula
\begin{equation}\label{e:DN0}
\lim_{R\to\infty} {R}^{h_{Y}}\cdot \zd\Rr_{R}
=2^{\z_{\Delta_Y}(0)} \zd^* \sqrt{\Delta_Y}\cdot
\det\big(\frac{\Id-C_{12}}{2}\big) \ \, .
\end{equation}
\end{thm}

\bigskip

Now we can use Theorem \ref{t:main thm}, the BFK formula
\eqref{e:BFK0} and Theorem \ref{t:DN0} to obtain the local
invariant $C(Y)$ as a byproduct of our main theorems.

\bigskip

\begin{cor}\label{c:c9y}
The BFK constant $C(Y)$ in the case of Dirichlet boundary
condition is equal to
\begin{equation}\label{e:c(y)}
C(Y) = 2^{-\z_{\Del_Y}(0) - h_Y} \ \, .
\end{equation}
\end{cor}

\bigskip

This result is also proved in \cite{Lee} independently using the
local computation of symbols of $\Rr_R$.

\bigskip

In Section 5 we discuss the proof of the technical result which
was used in Section 4 in the computation of the adiabatic limit of
the $\z$-determinant of $\Rr_R$ . Our approach is based on the
representation of the inverse of $\Delta_R$ in terms of the heat
kernel $e^{-t\Delta_R}$ , which enables us to apply the heat
kernel analysis and some results proved in the first part of the
paper.
\bigskip

{\bf Acknowledgment } The first author wishes to express his
gratitude to Werner M\"uller for helpful discussions. The authors
also thank the referee for corrections and helpful suggestions,
all of which considerably improved this paper. A part of this work
was done during the first author's stay at MPI. He also wishes to
express his gratitude to MPI for financial support and various
help.

\vskip 1cm

\section{Small eigenvalues and scattering matrices}
\label{s:scat}

\bigskip

In this section we study the relation between the \emph{s-values}
of the operators $\Del_R$, $\Del_{i,R}$ and the scattering
matrices $C_i(\lambda)$ determined by the operators
$\Del_{i,\infty}$ on $M_{i,\infty}$. This analysis is necessary in
order to determine the large time contribution in the adiabatic
decomposition formula. The corresponding result for Dirac
Laplacians was formulated and proved in \cite{JPKW5}. Here we
treat the case of a general Laplace type operator and we need to
rework some of the details of the analysis presented in
\cite{JPKW5}.

\bigskip
 Now let $\psi$ be an element of
$\ker(\Del_Y)$ and $\la$ denote a sufficiently small real number.
The couple $(\psi , \la)$ determines a generalized eigensection
$E(\psi, \lambda) \in C^{\infty}(M_{1,\infty}, E)$ of the operator
$\Delta_{1,\infty}$ such that
$$
\Del_{1,\infty} E(\psi, \lambda)= \sql E(\psi, \lambda) \ \ .
$$
The function $\la \to E(\psi,\la)$ has a meromorphic extension to
a certain subset of $\bC$, in particular, this function is
analytic function on the interval $(-\delta, \delta)$ for
sufficiently small $\delta >0 $. The generalized eigensection
$E(\psi,\la)$ has the following expression on the cylinder
$[0,\infty)_u\times Y$,
\begin{equation}\label{e:ge1}
E(\psi, \la)=e^{-i\la u}\psi + e^{i\la u}C_1(\lambda)\psi
+\hat{E}({\psi},\la) \ \, ,
\end{equation}
where $\hat{E}({\psi},\la)$ is a smooth $L^2$-section orthogonal
to $\ker(\Del_Y)$ and $\hat{E}(\psi,\la)|_{u=R}$ and $\partial_u
\hat{E}(\psi,\la)|_{u=R}$ are exponentially decaying as
$R\to\infty$. The scattering matrix
$$C_1(\la) : \ker(\Delta_Y) \to \ker(\Delta_Y)$$
is a unitary operator. The analyticity of $E(\psi,\la)$ implies
that $\{C_1(\la)\}_{\la\in(-\delta,\delta)}$ is an analytic family
of linear operators. The operator $C_1(\la)$ satisfies the
following functional equation
\begin{equation}\label{e:ge2}
C_1(\la)\ C_1(-\la)=\Id \ \, .
\end{equation}
In particular, $C_1(0)^2=\Id$,  hence $C_1(0)$ is an involution
over $\ker(\Delta_Y)$.

\bigskip

 Let $\Phi_R$ be a normalized
eigensection of $\Del_{1,R}$ for the Dirichlet boundary problem,
which corresponds to the \emph{s-value} $\sql = \la(R)^2$ with
$|\la| \le R^{-\kappa}$ for some fixed $\kappa$ with $0 < \kappa
\le 1$. That is,
\begin{equation}\label{e:dbc}
\Del_{1,R}\Phi_R=\sql\Phi_R  \ \ , \ \
 \Phi_R|_{\{R\}\times Y}=0 \qquad  \text{and}\qquad
||\Phi_R||=1\ \ .
\end{equation}
The section $\Phi_R$ can be represented in the following way on
$[0,R]_u\times Y \subset M_{1,R}$
$$
\Phi_R=e^{-i\la u}\psi_1 + e^{i\la u}\psi_2 +\hat{\Phi}_R
$$
where $\psi_i \in \ker(\Del_Y)$ and $\hat{\Phi}_R$ is orthogonal
to $\ker(\Del_Y)$.

\bigskip

We introduce $F:=\Phi_R - E(\psi_1,\lambda)|_{M_{1,R}}$ where
$\la$ is the positive square root of $\la^2$ . Green's theorem
gives
\begin{align}\label{e:green}
0 =&\, \langle \Del_{1,R}F, F\rangle_{M_{1,R}} -
\langle F, \Del_{1,R}F \rangle_{M_{1,R}} \\
=& -\int_{\partial{M_{1,R}}} \langle {\partial_u} F|_{u=R},
F|_{u=R} \rangle \ dy +  \ \ \int_{\partial{M_{1,R}}} \langle
F|_{u=R},{\partial_u} F|_{u=R} \rangle\  dy \ \ , \notag
\end{align}
and we can obtain the following equalities
\begin{align}\label{e:key1}
&\ \ \ 2\lambda\ i \ \|\, C_1(\la)\psi_1 - \psi_2\, \|^2 \\
=& - \ \langle\  {\partial_u}(\, \hat{\Phi}_R- \hat{E}(\psi_1,\la)
\, )|_{u=R}\ ,\
(\, \hat{\Phi}_R-\hat{E}(\psi_1,\la)\, )|_{u=R} \ \rangle \notag \\
&+ \  \langle\ (\, (\hat{\Phi}_R-\hat{E}(\psi_1,\la)\, )|_{u=R}\
,\ {\partial_u}(\, \hat{\Phi}_R-
\hat{E}(\psi_1,\la)\, )|_{u=R} \ \rangle  \notag \\
=& - \ \langle\  {\partial_u} (\, \hat{\Phi}_R
-\hat{E}(\psi_1,\la)\, )|_{u=R}\ ,\  -
\hat{E}(\psi_1,\la)|_{u=R} \ \rangle \notag \\
&+ \  \langle\ -\hat{E}(\psi_1,\la) |_{u=R}\ ,\ {\partial_u} (\,
\hat{\Phi}_R- \hat{E}(\psi_1,\la)\, )|_{u=R} \ \rangle  \ \ .
\notag
\end{align}
The following lemma will be used to show that the right side of
(\ref{e:key1}) is exponentially small as $R\to\infty$,

\bigskip

\begin{lem}\label{l:est-dir}
For $R \gg 0 $, there exists a constant $C$ independent of $R$
such that
\begin{equation*}
||\, \partial_u \hat{\Phi}_R |_{u=R}\, ||_{Y} \le C .
\end{equation*}
\end{lem}

\begin{proof}
We have the representation of $\hat{\Phi}_R$ on the cylinder $[0,
R ]_u\times Y \subset M_{1,R}$,
$$
\hat{\Phi}_R(u,y)\ = \ \sum_{k=h_Y+1}^{\infty} (a_k(R)
e^{\sqrt{\mu_k^2-\la^2} u}+b_k(R)
e^{-\sqrt{\mu_k^2-\la^2}u})\phi_k
$$
where $\{\mu_k^2,\phi_k\}$ is the spectral resolution of
$\Delta_Y$ , such that $\{\phi_k\}_{k=1}^{h_Y}$ is an orthonormal
basis of $\ker(\Del_Y)$ . The normalized condition for $\Phi_R$
implies the inequality
$$
\sum_{k=h_Y+1}^{\infty}\int^R_{0} |\ a_k(R)e^{\sqrt{\mu_k^2-\la^2}
u}+ b_k(R) e^{-\sqrt{\mu_k^2-\la^2} u}\ |^2 \ du \ \leq 1 \ \,
$$
which leads to
\begin{align*}
1 \ \geq  & \ \sum_{k=h_Y+1}^{\infty} \biggl(\
\frac{1}{2\sqrt{\mu_k^2-\la^2}} \big(
\ |a_k(R)|^2(e^{2\sqrt{\mu_k^2-\la^2}R}-1) \\
& \qquad \qquad \qquad \qquad + |b_k(R)|^2(1 -
e^{-2\sqrt{\mu_k^2-\la^2} R})\ \big) + \ 2\Re (a_k(R){b_k}(R))R\
\biggr) \ \, .
\end{align*}
The boundary condition put the following constraint on the
coefficients $a_k(R)$, $b_k(R)$
\begin{equation*}
a_k(R)e^{\sqrt{\mu_k^2-\la^2}R}+ b_k(R)
e^{-\sqrt{\mu_k^2-\la^2}R}\ = \ 0 \ \ .
\end{equation*}
As a result, if $R\gg 0$, the following estimate holds,
\begin{align}\label{e:est0}
1 \ &\geq \ \sum_{k=h_Y+1}^{\infty} \ \frac{ |a_k(R)|^2
e^{2\sqrt{\mu_k^2-\la^2}R} } {4\sqrt{\mu_k^2-\la^2}} \ \biggl( \
1+  e^{2\sqrt{\mu_k^2-\la^2} R}
- 8 \sqrt{\mu_k^2-\la^2}R \ \biggr) \\
&\geq \ \sum_{k=h_Y+1}^{\infty} \ { (\mu^2_k-\la^2)|a_k(R)|^2} \
e^{2\sqrt{\mu_k^2-\la^2}R} \biggl(\ \frac{1+
e^{\sqrt{\mu_k^2-\la^2} R}
}{4(\mu_k^2-\la^2)^{3/2}} \ \biggr) \notag \\
&\geq \ \sum_{k=h_Y+1}^{\infty} \ { (\mu^2_k-\la^2)|a_k(R)|^2} \
e^{2\sqrt{\mu_k^2-\la^2}R} \ \, . \notag
\end{align}
On the other hand, we can see that
\begin{align}\label{e:est1}
||\ \partial_u \hat{\Phi}_R |_{u=R}\ ||_{Y}^2 \ = \ 4
\sum_{k=h_Y+1}^{\infty} ({\mu_k^2-\la^2}) |a_k(R)|^2
e^{2\sqrt{\mu_k^2-\la^2} R} \ \ .
\end{align}
By (\ref{e:est0}) and (\ref{e:est1}), there is a constant $C$
independent of $R$ such that
$$
\|\ \partial_u \hat{\Phi}_R |_{u=R}\ \|_{Y} \le C \ \ .
$$

\end{proof}

\bigskip

Now Lemma \ref{l:est-dir} and the fact that
$\hat{E}(\psi,\la)|_{u=R}$ and $\partial_u
\hat{E}(\psi,\la)|_{u=R}$ are exponentially decaying as
$R\to\infty$ imply
\begin{equation}\label{e:scat1}
\|\, C_1(\la)\psi_1 - \psi_2\, \|^2 \le c_1\lambda^{-1} e^{-c_2R}
\le e^{-c_3R} \ \,
\end{equation}
for some positive constants $c_1, c_2, c_3$ . The second
inequality follows from the condition \eqref{Condition A}. Now the
Dirichlet boundary condition at $u=R$ of
$$
\Phi_R=e^{-i\la u}\psi_1 + e^{i\la u}\psi_2 +\hat{\Phi}_R
$$
provides us with the following equality,
\[
\psi_2= -e^{-2i\la R} \psi_1.
\]
From this equality and the estimate (\ref{e:scat1}), we get the
following inequality,
\begin{align}\label{e:sc1}
&\|\, e^{2i\la R}C_1(\la)\psi_1 + \psi_1\, \| \le e^{-cR} \ \ .
\end{align}
Recall that $\{ C_{1}(\lambda)\}_{\lambda\in(-\delta,\delta)}$ is
an analytic family of the operators. Analytic perturbation theory
guarantees the existence of the real analytic functions
$\alpha_j(\lambda)$ of $\lambda\in(-\delta,\delta)$ , such that
$\exp(i\alpha_j(\lambda))$ are the corresponding eigenvalues of
$C_{1}(\lambda)$ for $\lambda \in (-\delta, \delta)$ . Hence, from
(\ref{e:sc1}), we can obtain
\begin{align*}\label{e:sca1}
|\, e^{i(2\la R+\alpha_j(\la))} + 1\, | \le  e^{-cR} \ \ .
\end{align*}
This immediately implies

\bigskip

\begin{prop}\label{p:small1} For $R\gg 0$,
the positive square root $\la(R)$ of \emph{s-value} $\lambda(R)^2$
of $\Del_{1,R}$ with $\lambda(R)\le R^{-\kappa}$ $(0<\kappa\le 1)$
satisfies
\begin{equation}\label{e:2rl-p}
2R\lambda(R)+\alpha_j(\lambda(R))= (2k+1)\pi +O(e^{-cR}) \ \
\end{equation}
for an integer $k$ with $0 < (2k+1)\pi - \alpha_j(\lambda(R))\le
R^{1-\kappa}$, where $\exp(i\alpha_j(\lambda))$ is an eigenvalue
of the  unitary operator $C_1(\lambda): \ker(\Del_Y) \to
\ker(\Del_Y)$ .
\end{prop}

\bigskip

Now, we consider equation (\ref{e:2rl-p}) when $k=0$. The function
$\alpha_j(\la)$ is a real analytic function of $\la$ , hence we
have
\begin{equation}\label{e:exclude}
2R\lambda(R)+\alpha_{j0}+\alpha_{j1}\la(R) +\alpha_{j2}\la(R)^2
+\cdots = \pi + O(e^{-cR}) \ \,
\end{equation}
for some constants $\alpha_{jk}$'s. The operator $C_1(0)$ is an
involution, so $\alpha_{j0}=0$ or $\alpha_{j0}=\pi$ . It is not
difficult to show that, if we assume $\alpha_{j0}=\pi$ , then
$\la$ decays exponentially. However, the operator $\Del_{1,R}$
does not have the exponentially decaying eigenvalues, therefore
$\alpha_{j0}=0$. Now we proved

\bigskip

\begin{prop}\label{p:model0} For $R\gg 0$,
the positive square root $\la(R)$ of $s-value$ $\lambda(R)^2$ of
$\Del_{1,R}$ with $\la(R)\le R^{-\kappa}$ $(0<\kappa\le 1)$
satisfies
\begin{equation}\label{e:model1}
2R\lambda(R)= (2k+1)\pi + O({R^{-\kappa}}) \qquad \text{or} \qquad
2R\lambda(R)= 2k\pi + O( {R^{-\kappa}})
\end{equation}
where  $0 < (2k+1)\pi \le R^{1-\kappa}$ or $\ 0 <2k\pi \le
R^{1-\kappa}$ .
\end{prop}

\bigskip
Now one can easily prove that the similar result as in Proposition
\ref{p:model0} holds for $\Del_{2,R}$ simply repeating the
previous argument with the scattering matrix
$C_2(\la):\ker(\Del_Y)\to \ker(\Del_Y)$.

\bigskip

We are going to formulate Proposition \ref{p:model0} and the
corresponding result for $\Del_{2,R}$ in terms of certain model
operator over $S^1$. Let $U : W \to W$ denote a unitary operator
acting on a $d$-dimensional vector space $W$ with eigenvalues
$e^{i\alpha_j}$ for $j=1,\cdots, d$ . We define the operator
$\Del(U)$,
$$
\Del(U):=-\frac1 4\frac{d^2}{du^2} : C^{\infty}(S^1 , E_U) \to
C^{\infty}(S^1 , E_U)
$$
where $E_U$ is the flat vector bundle over
$S^1=\mathbb{R}/\mathbb{Z}$ defined by the holonomy ${U}$ . The
spectrum of $\Del(U)$ is equal to
\begin{equation}\label{e:spec}
\{ \ (\pi k+\frac{1}{2}\alpha_j)^2 \ | \ k\in\mathbb{Z} \ , \
j=1,\cdots, d \ \} \ \, .
\end{equation}
We also have
\begin{equation}\label{e:s1}
\zd\Del(U)=4^d\prod^d_{j=1} \sin^2(\frac{\alpha_j}{2}) \ \,
\end{equation}
if $\alpha_j\neq 2k\pi$ ($k\in\mathbb{Z}$) for $j=1,\cdots, d$
(see for instance \cite{LP2}). Putting $\overline{C}_i:=-C_i(0)$,
by definition, the operator $\Del(\overline{C}_i)$ has a
nontrivial kernel which is determined by $(1)$-eigenspace of
$\overline{C}_i$. We denote by $h_i$ the dimension of this space.

\bigskip

\begin{prop}\label{p:model1}
For any family of eigenvalues $\lambda(R)^2$ of $\Del_{i,R}$
converging to zero as $R \to \infty$, there exists the eigenvalue
$\lambda^2_k$ of $\Del(\overline{C}_i)$ with $\la_k>0$ so that for
$R\gg 0$,
\begin{equation}\label{e:m1}
R^2\lambda(R)^2= \la_k^2  + O({R^{1-2\kappa}}) \ \, ,
\end{equation}
and there is $R_1$ depending on $R$ with
$|R_1^{1-\kappa}-R^{1-\kappa}|\le \frac{\pi}{2}$ such that
\eqref{e:m1} defines one to one correspondence between the
eigenvalues of $\Del_{i,R}$ with $0<\la(R)^2\le R^{-2\kappa}$ and
the eigenvalues of $\Del(\overline{C}_i)$ with $0<\la_k^2\le
R_1^{2-2\kappa}$ and $\la_k>0$.
\end{prop}

\begin{proof}The equality \eqref{e:m1} follows
from Proposition \ref{p:model0}, the corresponding result for
$\Del_{2,R}$ and the definition of $\Del(\overline{C}_i)$. For the
second statement, by definitions, it is obvious that \eqref{e:m1}
defines an injective map from the eigenvalues of $\Del_{i,R}$ with
$0<\la(R)^2\le R^{-2\kappa}$ to the eigenvalues of
$\Del(\overline{C}_i)$ with $0<\la_k^2\le R^{2-2\kappa}$ and
$\la_k>0$. To define $R_1$ with the desired property, let us
decompose $M_{i,R}$ into $M_i$ and the cylindrical part of length
$R$. Then the restrictions of $\Del_{i,R}$ onto these decomposed
parts  provide us with the Laplace type operators imposing the
Dirichlet boundary conditions. By the mini-max principle, for
$R\gg 0$, the number of eigenvalues $\le R^{-2\kappa}$ of
$\Del_{i,R}$ is same as the number of eigenvalues $\le
R^{-2\kappa}$ of the operator over the cylindrical part since
there are no such small eigenvalues of the operator over $M_i$. By
the explicit computation over the cylinder of length $R$, the
eigenvalues of the operator over the cylinder of length $R$ are
given by $h_Y$-copies of ${k^2\pi^2}R^{-2}$ with $k\in\mathbb{N}$.
Therefore, the number of eigenvalues $\le R^{-2\kappa}$ of the
operator over the cylindrical part is given by
$h_Y[\pi^{-1}R^{1-\kappa}]$. Using \eqref{e:spec}, we can choose
$R_1$ such that $|R_1^{1-\kappa}-R^{1-\kappa}|\le \frac{\pi}{2}$
and  $h_Y[\pi^{-1}R^{1-\kappa}]$ is same as the number of the
eigenvalues of $\Del(\overline{C}_i)$ with $\la_k^2\le
R_1^{2-2\kappa}$ and $\la_k>0$. This completes the proof.
\end{proof}

\bigskip

Now we split
\[
\Tr(e^{-tR^2\Del_{i,R}})= \Tr_{1,R}(e^{-tR^2\Del_{i,R}}) +
\Tr_{2,R}(e^{-tR^2\Del_{i,R}}) \ \ ,
\]
where $\Tr_{1,R}(\cdot)$, $\Tr_{2,R}(\cdot)$ denote the parts of
the traces restricted to the nonzero eigenvalues $> R^{\frac12}$
or $\le R^{\frac12}$ of $R^2\Del_{i,R}$ respectively. Similarly,
we split
\[
\Tr(e^{-t\Del(\overline{C}_{i})}) -h_i
=\Tr_{1,R}(e^{-t\Del(\overline{C}_{i})}) +\Tr_{2,
R}(e^{-t\Del(\overline{C}_{i})}) \] where $\Tr_{1,R}(\cdot)$,
$\Tr_{2,R}(\cdot)$ denote the parts of the traces restricted to
the nonzero eigenvalues $> R_1^{\frac12}$ or $\le R_1^{\frac12}$
of $\Del(\overline{C}_{i})$ respectively. Now we have the estimate
for $\Tr_{2,R}(\cdot)$ in the following proposition.

\bigskip

\begin{prop}\label{p:diff}
For $R\gg 0$ , there exist positive constants $c_1,c_2$ such that

$$|\ \Tr_{2,R}(e^{-tR^2\Del_{i,R}})
-\frac12[\, \Tr_{2,R}(e^{-t\Del(\overline{C}_{i})}) - h_i \, ]\ |
\ \le\ c_{1}\,{R^{-\frac 14}}\, t\, e^{-c_{2} t} \ \ .$$
\end{prop}

\begin{proof}
We apply Proposition \ref{p:model1} for fixed $\kappa=\frac 34$
and obtain that for any eigenvalue $\lambda(R)^2$ of $\Del_{i,R}$
with $|\lambda(R)|\le R^{-\frac 34}$, there exists a function
$\alpha(R)$ such that
$$
R^2\lambda(R)^2= \lambda_j^2 + \alpha(R), \qquad |\alpha(R)| \le
c\, R^{-\frac 12}$$ if $R$ is sufficiently large. We use the
elementary inequality $|e^{-\lambda}-1|\le |\lambda|e^{|\lambda|}$
to get
\begin{multline*}
|e^{-tR^2\lambda(R)^2}-e^{-t\lambda_j^2}|= |e^{-t\lambda_j^2}
(e^{-t[R^2\lambda(R)^2-\lambda_j^2]}-1)|\\
\le c\, R^{-\frac12} \, t\, e^{-(\lambda_j^2-\alpha(R))t} \le c\,
{R^{-\frac 12}}\, t\,  e^{-\frac 12 \lambda_j^2t} \ \ .
\end{multline*}
Let us fix a sufficiently large $R$. We take the sum over finitely
many nonzero eigenvalues $\lambda(R)^2$ of $\Del_{i,R}$ with
$\lambda(R)^2\le R^{-\frac 32}$, and obtain
$$
|\ \Tr_{2,R}(e^{-tR^2\Del_{i,R}}) -\frac12[\,
\Tr_{2,R}(e^{-t\Del(\overline{C}_{i})}) -h_i \, ] \ |\le c\,
{R^{-\frac 12}}\, t\, \sum_{\lambda_j^2\le R_1^{\frac 12}}
e^{-\frac 12\lambda_j^2t} \ \ .$$ The operator
$\Del(\overline{C}_{i})$ is a Laplace type operator over $S^1$ ,
hence the number of eigenvalues $\lambda_j^2$ with $\lambda_j^2\le
R_1^{\frac 12}$ can be estimated by $R_1^{\frac 14}$. Since
$|R_1^{\frac14}-R^{\frac14}|\le \frac{\pi}{2}$, we have
$$
c \, {R^{-\frac 12}}\,t\, \sum_{\lambda_j^2 \le R_1^{\frac 12}}
e^{-\frac 12\lambda_j^2t} \le c'\, {R^{-\frac 14}}\, t\,
e^{-\frac12\lambda^2_1 t}
$$
where $\la_1^2$ denotes the first non-zero eigenvalue of
$\Del(\overline{C}_{i})$. This completes the proof.
\end{proof}

\bigskip

Now we shall prove the corresponding result for the
\emph{s-values} of $\Del_R$ over $M_R$. Let $\Psi_R$ denote (a
normalized) eigensection of $\Del_R$ corresponding to
\emph{s-value} $\la^2$, that is, $\Del_R\Psi_R =\la^2\Psi_R$ and
$\|\Psi_R\| = 1$. Over the cylindrical part $[-R,R]_u\times Y$ in
$M_R$, the eigensection $\Psi_R$ corresponding to \emph{s-value}
$\la^2$ of $\Del_R$ has the following form,
\begin{equation}\label{e:psi}
\Psi_R = e^{-i\la u}\psi_1 + e^{i\la u}\psi_2 +\hat{\Psi}_R
\end{equation}
where $\psi_i\in \ker(\Del_Y)$ and $\hat{\Psi}_R$ is orthogonal to
$\ker(\Del_Y)$. We first need the following lemma, where
$\{0\}\times Y$ denotes the cutting hypersurface in $M_R$.

\bigskip

\begin{lem} \label{l:cut-est}
We have the following estimates
\begin{align*}
||\hat{\Psi}_R|_{u=0}||_Y \leq  c_1 e^{-c_2R} \quad, \quad
||\partial_u\hat{\Psi}_R|_{u=0}||_Y \leq  c_1e^{-c_2R}
\end{align*}
where $c_1,c_2$ are positive constants independent of $R$ .
\end{lem}

\begin{proof}
The section $\hat{\Psi}_R$ has the following form on
$[-R,R]_u\times Y \subset M_R$,
\begin{align*}
\hat{\Psi}_R(u,y)\ = \ \sum^{\infty}_{k=h_Y+1} (a_k(R)
e^{\sqrt{\mu_k^2-\la^2} u}+ b_k(R)
e^{-\sqrt{\mu_k^2-\la^2}u})\phi_k \ \ .
\end{align*}
The normalization condition on the eigensection implies
$$
\sum_{k=h_Y+1}^{\infty}\int^R_{-R} |\
a_k(R)e^{\sqrt{\mu_k^2-\la^2} u}+ b_k(R) e^{-\sqrt{\mu_k^2-\la^2}
u}\ |^2 \ du \ \leq 1 \ \, ,
$$
and now we have the following estimates for sufficiently large $R$
\begin{align*}
1 \ \geq  & \ \sum_{k=h_Y+1}^{\infty} \biggl(\
\frac{1}{2\sqrt{\mu_k^2-\la^2}}\ [ \
|a_k(R)|^2(e^{2\sqrt{\mu_k^2-\la^2}R}-
e^{-2\sqrt{\mu_k^2-\la^2}R} ) \\
& \quad \quad + |b_k(R)|^2(e^{2\sqrt{\mu_k^2-\la^2}R} -
e^{-2\sqrt{\mu_k^2-\la^2} R})\ ] +
\ 4\Re (a_k(R){b_k}(R))R\  \biggr) \\
\geq &  \ \sum_{k=h_Y+1}^{\infty} \
\frac{1}{4\sqrt{\mu_k^2-\la^2}} \ ( \
\ |a_k(R)|^2 e^{2\sqrt{\mu_k^2-\la^2}R} \\
& \quad \quad \qquad + |b_k(R)|^2
e^{2\sqrt{\mu_k^2-\la^2}R}  - 16|a_k(R){b_k}(R)| R \ ) \\
\geq &  \ \sum_{k=h_Y+1}^{\infty} \
\frac{1}{8\sqrt{\mu_k^2-\la^2}}\ ( \ \ |a_k(R)|^2
e^{2\sqrt{\mu_k^2-\la^2}R} + |b_k(R)|^2 e^{2\sqrt{\mu_k^2-\la^2}R}
\ ) \ \ .
\end{align*}

This immediately implies
\begin{align*}
\sum_{k=h_Y+1}^{\infty} |a_k(R)|^2  + |b_k(R)|^2 \ \leq \
c_1e^{-\sqrt{\mu_{h_Y+1}^2 -\la^2}R} \ \leq \ c_1e^{-c_2 R}
\end{align*}
for some positive constants $c_1,c_2$. Hence, the first estimate
is proved and the proof of the second estimate follows in the same
way.
\end{proof}

\bigskip

 Changing variable $v=u+R$, we
regard that the cylindrical part is given by $[0,2R]_v\times Y$.
In particular, we have the new expression for $\Psi_R$ from
\eqref{e:psi},
\[
\Psi_R = e^{-i\la v}\phi_1^1 + e^{i\la v}\phi_2^1 +\hat{\Psi}_R
\]
where $\phi_1^1=e^{i\la R}\psi_1$, $\phi_2^1=e^{-i\la R}\psi_2$.
Now repeating the argument which leads us to \eqref{e:scat1}, we
obtain
\begin{equation}\label{c1}
\|\, C_1(\la)\phi_1^1 - \phi_2^1\, \|  \le e^{-cR} \ \,
\end{equation}
for a positive constant $c$. Note that here we used the condition
\eqref{Condition A} and Lemma \ref{l:cut-est}. Now we want to get
the corresponding estimate involving the scattering matrix
$C_2(\la)$. For this, we change the variable by $v=u-R$ and regard
the cylindrical part as $[-2R,0]_v\times Y$. Then we have the
corresponding expression for $\Psi_R$,
\[
\Psi_R = e^{-i\la v}\phi_1^2 + e^{i\la v}\phi_2^2 +\hat{\Psi}_R
\]
where $\phi_1^2=e^{-i\la R}\psi_1$, $\phi_2^2=e^{i\la R}\psi_2$.
We again repeat the previous argument to obtain
\begin{equation}\label{c2}
\|\, C_2(\la)\phi_2^2 - \phi_1^2\, \|  \le e^{-cR} \ \,
\end{equation}
for a positive constant $c$. Here $C_2(\la)$ is the scattering
matrix defined from the generalized eigensection attached to
$(\la, \phi^2_2)$. By definition, we have
\begin{equation}\label{e:phi}
\phi^1_1=e^{2i\la R}\phi^2_1\ \ , \quad  \phi^1_2=e^{-2i\la
R}\phi^2_2\ \ .
\end{equation}
Now, combining \eqref{c1}, \eqref{c2} and \eqref{e:phi}, we get
\begin{equation}\label{c12}
\|\, e^{4i\la R} C_1(\la)\circ C_2(\la)\phi_2^1 - \phi_2^1\, \|
\le e^{-cR} \ \ .
\end{equation}
As before, $C_1(\la)\circ C_2(\la)$ is an analytic family for
$\la\in(-\delta,\delta)$ for sufficiently small $\delta>0$. Then
there exist the analytic functions $\alpha_j(\la)$ for
$\la\in(-\delta,\delta)$ such that $\exp((i\alpha_j(\lambda))$ are
the eigenvalues of the  unitary operator $
C_{12}(\la):=C_1(\la)\circ C_2(\la)$ on $\ker(\Del_Y)$. Hence the
equality \eqref{c12} implies
\[
|\, e^{i(4\la R+\alpha_j(\la))} - 1\, | \le e^{-cR} \ \ .
\]
Therefore we obtain

\bigskip

\begin{prop}\label{p:small eigen1}
For $R\gg 0$, the positive square root $\la(R)$ of \emph{s-value}
$\lambda(R)^2$ of $\Del_R$ with $\la(R)\le R^{-\kappa}$ satisfies
\begin{equation}\label{e: s1}
4R\lambda(R) + \alpha_j(\lambda(R)) =2k\pi + O(e^{-cR})
\end{equation}
for an integer $k$ with $0 < 2k\pi-\alpha_j(\la(R)) \le
4R^{1-\kappa}$ , where $\exp(i\alpha_j(\lambda))$ is the
eigenvalue of the  unitary operator $C_{12}(\la)$ on
$\ker(\Del_Y)$.
\end{prop}

\bigskip

\begin{rem}\label{r:c12}
Note that the spectrum of the unitary operator $C_{12}:=C_{12}(0)$
acting on $\ker(\Del_Y)$ consists of $m$ eigenvalues of $-1$ (
such that $h_Y - m \geq 0$ is an even number ) and $\{ \
e^{i\alpha_j(0)}, e^{-i\alpha_j(0)} \ | \ j=1, \cdots, \frac{h_Y -
m}{2} \ \}$ where $\alpha_j(0)$ is not equal to $k\pi$ for
$k\in\mathbb{Z}$ . This follows from the argument presented around
\eqref{e:exclude} and the condition \eqref{Condition A}.
\end{rem}

\bigskip
Now we follow the way to prove Proposition \ref{p:model1} and
obtain
\bigskip

\begin{prop}\label{p:model2}
For any family of eigenvalues $\lambda(R)^2$ of $\Del_{R}$
converging to zero as $R \to \infty$, there exists the eigenvalue
$\lambda^2_k$ of $\Del(C_{12})$ with $\la_k>0$ so that for $R\gg
0$,
\begin{equation}\label{e:m2}
4R^2\lambda(R)^2= \la_k^2  + O({R^{1-2\kappa}}) \ \, ,
\end{equation}
and there is $R_1$ depending on $R$ with
$|R_1^{1-\kappa}-R^{1-\kappa}|\le \frac{\pi}{4}$ such that
\eqref{e:m2} defines one to one correspondence between the
eigenvalues of $\Del_{R}$ with $0<\la(R)^2\le R^{-2\kappa}$ and
the eigenvalues of $\Del(C_{12})$ with $0<\la_k^2\le
4R_1^{2-2\kappa}$ and $\la_k>0$.
\end{prop}

\bigskip

We split
\[
\Tr(e^{-tR^2\Del_{R}})= \Tr_{1,R}(e^{-tR^2\Del_{R}}) +
\Tr_{2,R}(e^{-tR^2\Del_{R}}) \ \ ,
\]
where $\Tr_{1,R}(\cdot)$, $\Tr_{2,R}(\cdot)$ denote the parts of
the traces restricted to the nonzero eigenvalues $>R^{\frac12}$ or
$\le R^{\frac12}$ of $R^2\Del_R$ respectively. Similarly we split
\[
\Tr(e^{-t\frac14\Del(C_{12})})=\Tr_{1,R}(e^{-t\frac14\Del(C_{12})})
+\Tr_{2, R}(e^{-t\frac14\Del(C_{12})}) \] where
$\Tr_{1,R}(\cdot)$, $\Tr_{2,R}(\cdot)$ denote the parts of the
traces restricted to the nonzero eigenvalues $>R_1^{\frac12}$ or
$\le R_1^{\frac12}$ of  $\frac14\Del(C_{12})$ respectively. As in
Proposition \ref{p:diff}, we can prove the following proposition.

\bigskip

\begin{prop}\label{p:diff12}
For $R\gg 0$ , there exist positive constants $c_1,c_2$ such that

$$|\ \Tr_{2,R}(e^{-tR^2\Del_{R}})
-\frac12 \Tr_{2,R}(e^{-t\frac14\Del(C_{12})}) \ | \ \le \
c_{1}\,{R^{-\frac 14}}\, t\, e^{-c_{2} t} \ \ .$$
\end{prop}

\vskip 1cm

\section{Proof of Theorem \ref{t:main thm}} \label{s:stc}

\bigskip

In this section we present a proof of Theorem \ref{t:main thm}.
Since the analysis of \emph{s-values} is done in Section 2, now we
can proceed by a standard way as in \cite{JPKW4} and \cite{JPKW5}.

\bigskip

We define relative $\z$-function $\zeta_{\text{rel}}^R(s)$,
\begin{equation}\label{e:rel1}
\zeta_{\text{rel}}^R(s):=\frac{1}{\Gamma(s)}\int^{\infty}_{0}t^{s-1}
\Tr(e^{-t\Del_R} - e^{-t\Del_{1,R}} - e^{-t\Del_{2,R}} ) \ dt \ \,
,
\end{equation}
and we decompose $\zeta_{\text{rel}}^R(s)$ into two parts
$$\z_{\text{s}}^R(s)=\frac{1}{\Gamma(s)}\int^{R^{2-\e}}_0
(\cdot)\ dt \ \ , \ \ \z_\text{l}^R(s)=\frac{1}{\Gamma(s)}
\int^{\infty}_{R^{2-\e}} (\cdot)\ dt \ \ $$ where $\e>0$ is a
fixed sufficiently small number. The derivatives of
$\z_\text{s}^R(s)$ and $\z_\text{l}^R(s)$ at $s=0$ give the small
and large time contributions to our formula. First, we prove

\bigskip

\begin{lem}\label{l:asymp} There exist positive
constants $c_1$ and $c_2$ such that
\begin{equation*}
|\, \Tr(e^{-t\Del_R} - e^{-t\Del_{1,R}} - e^{-t\Del_{2,R}}) -
\frac{1}{2} \Tr(e^{-t\Del_Y})\, | \le  c_1e^{-c_2\frac {R^2}t} \ \
.
\end{equation*}
\end{lem}

\begin{proof}
By the standard application of Duhamel principle as in
\cite{JPKW2}, \cite{JPKW5}, the estimate of $\Tr(e^{-t\Del_R} -
e^{-t\Del_{1,R}} - e^{-t\Del_{2,R}})$ follows from the estimate of
the parametrices of the heat kernels $e^{-t\Del_R}$ ,
$e^{-t\Del_{i,R}}$. These parametrices are constructed from the
heat kernels on the closed manifold $M_R$ and heat kernels of the
boundary problems on the half infinite cylinders. The interior
contributions cancel each other out up to the error term of the
size $O(e^{-c\frac {R^2}t})$ for a positive constant $c$ and only
the boundary contribution is left. This boundary term is equal to
\begin{multline*}
\int_{-R}^{R}{\frac1{\sqrt{4{\pi}t}}}\Tr(e^{-t\Del_Y})\, du -
2\int_{0}^{R}{\frac1{\sqrt{4{\pi}t}}}\{1 -
e^{-{\frac{u^2}{t}}}\}\Tr(e^{-t\Del_Y})\, du \\
=2\int_{0}^{R}{\frac1{\sqrt{4{\pi}t}}} e^{-{\frac{u^2}{t}}}
\Tr(e^{-t\Del_Y})\, du
=\frac1{\sqrt{{\pi}}}\int_{0}^{\frac{R}{\sqrt{t}}}
e^{-v^2}\Tr(e^{-t\Del_Y})\, dv \\ = \frac 12 \Tr (e^{-t\Del_Y}) +
O(e^{-\frac {R^2}t}) \ \ . \end{multline*} This completes the
proof.
\end{proof}

\bigskip
Now we can determine the small time part in \eqref{e:rel1}.

\bigskip

\begin{prop} \label{p:sc}
We have
\begin{equation*}
\lim_{R\to\infty}[\, ({{\z}_\text{s}^R})'(0) -
\frac{h_Y}{2}(\gamma + (2-\epsilon)\log R)\, ] = \frac{1}{2}
\z_{\Delta_Y}'(0) \ \ ,
\end{equation*}
where
$$\z_{\Delta_Y}(s) = \frac 1{\G(s)}
\int_0^{\infty}t^{s-1} \big(\Tr( e^{-t\Del_Y}) - h_Y\big)\ dt \ \
.$$
\end{prop}

\begin{proof}
By Lemma \ref{l:asymp}, the function
\begin{equation*}\label{e:asymp1}
f_R(s) = \frac 1{\G(s)} \int_0^{R^{2-\e}} t^{s-1}
\big(\Tr(e^{-t\Del_R} - e^{-t\Del_{1,R}} - e^{-t\Del_{2,R}}) -
\frac {1}{2} \Tr(e^{-t\Del_Y})\big)\ dt \ \,
\end{equation*}
is a holomorphic function of $s$ on the whole complex plane.
Moreover, the following equalities hold
\begin{equation*}\label{e:asymp2}
\lim_{R \to \infty}f_R(0) = 0 \quad ,  \quad \lim_{R \to \infty}
\frac d{ds}f_R(s) \Bigr|_{s=0} = 0 \ \, .
\end{equation*}
Combining these facts with the following equality
\begin{align}\label{contr1}
 \frac{d}{ds}\Bigr|_{s=0}  \biggl( \frac{h_Y}{\Gamma(s)}
\int^{R^{2-\e}}_{0}t^{s-1} \ dt\biggr) =  h_Y(\gamma +
{(2-\e)}\log R) \ \ ,
\end{align}
completes the proof.
\end{proof}

\bigskip

To deal with the large time part, we need the following lemma.
\bigskip

\begin{lem}\label{l:int}
For $R\gg 0$, there exists a positive constant $c_1$ such that
\begin{align*}
\int^{\infty}_{R^{-\epsilon}} t^{-1}
\Tr_{1,R}(e^{-tR^2\Del_{i,R}})\, dt \le c_1\, e^{-R^{\frac 12-\e}}
\end{align*}
and the similar estimates hold for $\Tr_{1,R}(e^{-tR^2\Del_{R}})$,
$\Tr_{1,R}(e^{-t\Del(\overline{C}_i)})-h_i$ and
$\Tr_{1,R}(e^{-t\frac14\Del(C_{12})})$.
\end{lem}

\begin{proof} Let $\la_{k_0}^2(R)$
denote the smallest \emph{large} eigenvalue of $\Del_{i,R}$ such
that $\la_{k_0}^2(R) > R^{-\frac32}$. Then, if $R\gg 0$ we have
\begin{align*}
&\Tr_{1,R}(e^{-tR^2\Del_{i,R}} )= \sum_{\la_k^2 >
R^{-\frac32}}e^{-tR^2\la_k^2} = \sum_{\la_k^2 >
R^{-\frac32}} e^{-(tR^2-1)\la_k^2}e^{-\la_k^2} \\
&\quad \le e^{-(tR^2-1)\la_{k_0}^2}\sum_{\la_k^2 >
R^{-\frac32}}e^{-\la_k^2} \le e^{-(tR^2-1)\la_{k_0}^2} \Tr (
e^{-\Del_{i,R}}) \\
&\quad\quad \le c_2\, R\, e^{-(tR^2-1)R^{-\frac32}} \le c_3\, R\,
e^{-{R}^{\frac12}t}
\end{align*}
for positive constants $c_2, c_3$. We have used here the obvious
estimate
$$\Tr( e^{-\Delta_{i,R}}) \le c\, \text{vol}(M_{i,R})
\le c' R \ \ $$ for positive constants $c,c'$. Now we have
\begin{align*}
&\int^{\infty}_{R^{-\e}}{{t}^{-1}}\Tr_{1,R}(e^{-tR^2\Del_{i,R}})\,
dt
\le \int^{\infty}_{R^{-\e}}{{t}^{-1}}c_3\, R\, e^{-t{R^{\frac12}}}\, dt\\
&\qquad \le c_3 R \int_{R^{\frac 12 - \e}}^{\infty}e^{-v}dv \le
c_1\, e^{-R^{\frac 12 - \e}}  \ \ .
\end{align*}
This completes the proof of the first estimate and the other cases
can be proved in the same way.

\end{proof}

\bigskip

Now we can express the large time part in terms of the model
operators.
\bigskip

\begin{prop}\label{t:lc}
\begin{equation*}\label{e:lc}
\lim_{R\to\infty} \int^{\infty}_{R^{2-\e}} t^{-1}\Tr(e^{-t\Del_R}
- e^{-t\Del_{1,R}} - e^{-t\Del_{2,R}}) \ dt + \frac
{h_Y}{2}(\gamma - \e{\cdot}\log R)
\end{equation*}
$$
=\frac12 \frac{d}{ds}\biggm|_{s=0}
\frac{1}{\Gamma(s)}\int^{\infty}_0 t^{s-1}\ \Bigr(
\Tr(e^{-\frac{t}{4}\Del({C_{12}})} -e^{-t\Del(\overline{C}_1)}-
e^{-t\Del(\overline{C}_2)})+ {h_Y} \Bigr) \ dt  \ \ .
$$
\end{prop}

\begin{proof} First, let us observe that Remark \ref{r:c12} and
the relation $C_i(0)^2=\mathrm{Id}$ imply $h_Y=h_1+h_2$. Using
this and the change of variable $t\to R^{-2}t$, one can obtain
following equality from Proposition \ref{p:diff}, \ref{p:diff12}
and Lemma \ref{l:int},
\begin{multline*}
\lim_{R\to\infty} \Bigr( \int^{\infty}_{R^{2-\e}}
t^{-1}\Tr(e^{-t\Del_R} -
e^{-t\Del_{1,R}} - e^{-t\Del_{2,R}})\ dt\\
- \frac12 \frac{d}{ds}\biggm|_{s=0} \frac{1}{\Gamma(s)}
\int^{\infty}_{R^{-\e}} t^{s-1}\ [
\Tr(e^{-\frac{t}{4}\Del({C_{12}})} \\ -e^{-t\Del(\overline{C}_1)}\
- e^{-t\Del(\overline{C}_2)}) + {h_Y} ]\ dt \Bigr) =0  \ .
\end{multline*}
Note that near $t=0$,
\[
|\, \Tr( e^{-\frac{t}{4} \Del({C_{12}})}
-e^{-t\Del(\overline{C}_1)} -e^{-t\Del(\overline{C}_2)})\, | \le
c\sqrt{t}
\]
for a positive constant $c$. By this estimate, one can easily show
\begin{multline*}
\lim_{R \to \infty} \Bigr( h_Y(\gamma - \e{\cdot}\log\, R) -
\frac{d}{ds}\biggm|_{s=0} \frac{1}{\Gamma(s)}
\int^{R^{-\e}}_0t^{s-1}\ [\Tr( e^{-\frac{t}{4} \Del({C_{12}})}
\\\qquad \qquad -e^{-t\Del(\overline{C}_1)} -e^{-t\Del(\overline{C}_2)}) + h_Y]\
dt  \Bigr) = 0 \ \ .
\end{multline*}
These complete the proof.
\end{proof}

\bigskip

Proposition \ref{p:sc} and \ref{t:lc} combined together lead to
the following equality
\begin{align}\label{final}
&\ \lim_{R \to \infty}\Bigr((\z_\text{s}^R)'(0) - \frac{h_Y}{2}(\g
+ (2-\e){\cdot}\log R) + (\z_\text{l}^R)'(0) +\frac {h_Y}{2}(\g -
\e{\cdot}\log R)\Bigr)\\
& \qquad =\ \frac 12 \Bigr ( \ \z_{\Del_Y}'(0) + \z_{\frac 14
\Del({C_{12}})}'(0)- \z_{\Del(\overline{C}_1)}'(0)-
\z_{\Del(\overline{C}_2)}'(0) \Bigr) \ \ .\notag
\end{align}
\bigskip

Now the following proposition gives the exact value of the large
time contribution,
\bigskip

\begin{prop}\label{p:det1}
We have
$$ \zd \frac 14\Del({C_{12}})= 2^{2h_Y}
\det(\frac{\mathrm{Id}-C_{12}}{2})^2 \quad , \quad
\zd^*\Del(\overline{C}_i) = 2^{2h_Y}\ \ .
$$
\end{prop}

\begin{proof}
The first equality follows directly from (\ref{e:s1}). For the
second one, the zeta function of $\Del(\overline{C}_i)$ is given
by
$$
\z_{\Del(\overline{C}_i)}(s) =   h_i\,
2\pi^{-2s}\sum^{\infty}_{k=1} k^{-2s}+ (h_Y-h_i)\, 2\pi^{-2s}
\sum^{\infty}_{k=0}(k+ \frac 12)^{-2s} \ \
$$
where $h_i$ is the dimension of $(+1)$-eigenspace of
$\overline{C}_i$ . Then the derivative of
$\z_{\Del(\overline{C}_i)}(s)$ at $s=0$ is equal to $- h_Y\log 4$.
This completes the proof of the second one.
\end{proof}

\bigskip

Finally we obtain Theorem \ref{t:main thm} using the equality
\eqref{final} and Proposition \ref{p:det1}.

\vskip 1cm

\section{The adiabatic limit of  $\zd\Rr_R$}

\bigskip

In this section we study the behavior of $\zd \Rr_R $ when
$R\to\infty$.

\bigskip

Let us describe the construction of $\Rr_R$. It is defined as the
composition of the following maps
$$ C^{\infty}(Y,E|_Y) \overset{I_g}\longrightarrow
C^{\infty}(Y,E|_Y) \oplus C^{\infty}(Y,E|_Y)
\overset{\mathcal{\Kk_{R}}}\longrightarrow C^{\infty}({\overline
{M}_R},E)
$$
$$
\overset{\g_1}\longrightarrow C^{\infty}(Y, E|_Y)\oplus
C^{\infty}(Y,E|_Y) \overset{I_f}\longrightarrow C^{\infty}(Y,
E|_Y).
$$
Here $I_g(\phi):=(\phi,\phi)$ and $\Kk_{R}$ is the Poisson
operator of the operator $\Delta_{1,R} \sqcup \Delta_{2,R}$ over a
manifold $\overline{M}_R:= M_{1,R} \sqcup M_{2,R}$. For
$(\Phi_1,\Phi_2)$ where $\Phi_i$ is a section over $M_{i,R}$, the
map $\g_1$ is given by  $ \g_1(s):= (\partial_u\big|_{Y_1} \Phi_1,
\partial_{u}\big|_{Y_2} \Phi_2 )$ and $I_f(\phi,\psi):=\phi-\psi $ . It is well
known that the operator
$$\Rr_R:=I_f\, \g_1\,\Kk_{R}\, I_g \ : \ C^{\infty}(Y,E|_Y) \to
C^{\infty}(Y, E|_Y)
$$ is an elliptic,
nonnegative, pseudo-differential operator of order $1$ . By
definition, the operator $\Rr_R$ can be written as
$$
\Rr_R \ = \ \Nn_{1,R}+\Nn_{2,R} \ \
$$
where $\Nn_{i,R}$ is the Dirichlet to Neumann operator for
$\Delta_R|_{M_{i,R}}$.

\bigskip

A careful analysis of the small eigenvalues enables us to compute
the scattering contribution to the adiabatic limit of the
$\z$-determinant of $\Rr_R$. Let us recall that $\{\mu_k^2,
\phi_k\}_{k\in\mathbb{N}}$ denotes the spectral resolution of the
operator $\Delta_Y$ with $h_Y=\dim \ker(\Del_Y)$. The equality
(\ref{e:ge2}) implies
$$C_i(0)C_i'(0)=C_i'(0)C_i(0) \ \ ,$$
hence we may choose $\phi_k$ (for $1\leq k \leq h_Y$) so that
$\phi_k$ is a normalized eigensection for both operators $C_i(0)$
and $C_i'(0)$. Now, we have

\bigskip

\begin{prop}\label{p:prop-est2}
For any couple $(\phi_m , \phi_n)$ with $1 \le m,n \le h_Y$,
\begin{equation*}
\begin{split}
\langle \Nn_{i,R}\phi_m , \phi_n \rangle \ = \
\begin{cases} &
\frac{1}{R}(1-\frac{\alpha}{2R})^{-1} \quad \text{if}\quad m=n \ ,
\ \ C_i(0)\phi_m=-\phi_m \ \ , \\
& O(e^{-cR}) \quad \text {if} \quad m\neq n \quad\text{or}\quad
C_i(0)\phi_m=\phi_m \ \
\end{cases}
\end{split}
\end{equation*}
where $C_i'(0)\phi_n=i\alpha\phi_n$ , that is, $i\alpha$ is the
eigenvalue of $C_i'(0)$ and  $c$ is a positive constant .
\end{prop}

\begin{proof}
We present a proof for the case of $i=1$ . The case for $i=2$ can
be proved in the same way. Let $\Phi_R$ denote a solution of the
problem
$$\Delta_{M_{1,R}} \Phi_R=0 \quad \ \text{and} \quad \
\Phi_R|_{Y}=\phi_m \ \ ,$$ hence
\begin{equation}\label{e:dnr1}
\partial_u\Phi_R|_{u=R} = \Nn_{1,R}\phi_m
\ \, .
\end{equation}
To simplify notation in the proof we skip the indices $m$ in
$\phi_m$ and $R$ in $\Phi_R$ . Let us define
$$\Phi(\phi, \lambda ):=e^{-i\lambda R}\Phi \ \ ,$$
for a small positive $\lambda$. For such a $\lambda$ and
$\psi:=\phi_n \in \ker(\Delta_Y)$ , there exists the generalized
eigensection $E(\psi,\lambda)$ over $M_{1,\infty}$, which has the
following form on the cylinder $[0,\infty)_u\times Y \subset
M_{1,\infty}$,
$$
E(\psi,\lambda)=e^{-i\lambda u}\psi + e^{i\lambda
u}C_1(\lambda)\psi + \hat{E}(\psi,\la)
$$
where $\hat{E}(\psi, \la)$ is a $L^2$-section. We also define
$$
G =
G(\phi,\psi,\lambda):=E(\psi,\lambda)|_{M_{1,R}}-\Phi(\phi,\lambda)
\ \ .
$$
An auxiliary section $G(\phi,\psi,\lambda)$ has the following
properties
$$\Delta_{1,R}G(\phi,\psi,\lambda)=\lambda^2
E(\psi,\lambda) \ \ ,$$
$$G|_{u=R} = e^{-i\la R}\psi +
e^{i\la R}C_1(\la)\psi - e^{-i\la R}\phi + O(e^{-cR}) \ \ ,$$
$$\partial_{u} G|_{u=R}   =
- i\la e^{-i\la R}\psi + i\la e^{i\la R}C_1(\la)\psi - e^{-i\la
R}\Nn_{1,R}\phi + O(e^{-cR}) \ \ .$$ Green's formula for $G$ reads
as
\begin{multline}\label{e:above}
 \langle \Delta_{1,R}G,G \rangle_{M_{1,R}} - \langle
G,\Delta_{1,R} G \rangle_{M_{1,R}}\\
  = -\langle {\partial_u}G|_{\{R\}\times Y},
G|_{\{R\}\times Y} \rangle_{\{R\}\times Y} + \langle
G|_{\{R\}\times Y},{\partial_u}G|_{\{R\}\times Y}
\rangle_{\{R\}\times Y} \ \ .
\end{multline}
The equation (\ref{e:above}) can be rewritten as follows
\begin{equation}\label{e:above1}
\begin{split}
&\lambda^2(\, \langle\Phi,E\rangle_{M_{1,R}}-\langle
E,\Phi\rangle_{M_{1,R}})\\
=\quad & e^{-2i\lambda R}\langle \Nn_{1,R} \phi,
C_1(\lambda)\psi\rangle_Y -e^{2i\lambda R}\langle
C_1(\lambda)\psi, \Nn_{1,R} \phi\rangle_Y \\
+&\, i\lambda e^{-2i\lambda R} \langle
\phi,C_1(\lambda)\psi\rangle_Y +i\lambda e^{2i\lambda
R} \langle C_1(\lambda)\psi, \phi \rangle_Y \\
+&\, \langle \Nn_{1,R}\phi,\psi\rangle_Y - \langle
\psi,\Nn_{1,R}\phi\rangle_Y  -\langle \Nn_{1,R}\phi, \phi\rangle_Y
+
\langle \phi,\Nn_{1,R}\phi\rangle_Y\\
     -&\, i\lambda \langle \phi,\psi
\rangle_Y -\, i\lambda \langle \psi,\phi \rangle_Y+O(e^{-cR}) \ \
.
\end{split}
\end{equation}
We differentiate both sides of the equality (\ref{e:above1}) at
$\lambda=0$ and obtain
\begin{equation}\label{e:stram}
\begin{split}
-&\, 2i R(\, \langle \Nn_{1,R} \phi, C_1(0)\psi\rangle_Y
+ \langle C_1(0)\psi, \Nn_{1,R} \phi\rangle_Y)\\
+&\, \langle \Nn_{1,R}\phi,C_1'(0)\psi\rangle_Y - \langle
C_1'(0)\psi,\Nn_{1,R}\phi\rangle_Y \\
+&\, i(\, \langle \phi,C_1(0)\psi\rangle_Y + \langle C_1(0)\psi,
\phi \rangle_Y) - i \langle \phi,\psi \rangle_Y- i \langle
\psi,\phi \rangle_Y = O(e^{-cR}) \ \ .
\end{split}
\end{equation}
Proposition \ref{p:prop-est2} follows easily from (\ref{e:stram}).
Let us consider for instance the case of
$$\phi = \psi = \phi_n \in \ker(C_1(0)+ 1)
\subset \ker( \Delta_Y) \ \ .$$ Then, the equation (\ref{e:stram})
is now
$$(2iR-i\alpha)(\, \langle \Nn_{1,R} \phi,
\phi\rangle_Y + \langle  \phi, \Nn_{1,R}\phi\rangle_Y ) = 4i+
O(e^{-cR}) \ \ ,$$ and this gives the following formula,
\begin{equation}\label{e:stram1}
\langle \Nn_{1,R} \phi, \phi\rangle_Y + \langle  \phi,
\Nn_{1,R}\phi\rangle_Y = \frac 2R (1-\frac {\alpha}{2R})^{-1} +
O(e^{-cR}) \ \, .
\end{equation}

\end{proof}

\bigskip

Let us also observe the following fact, which is an immediate
corollary of Proposition \ref{p:prop-est2}.

\bigskip

\begin{cor}\label{c:adn}
We have
\begin{align*}
\langle \Rr_R \phi , \phi \rangle  \ = \ O(e^{-cR})\quad
\text{for}\quad \phi \in \ker(C_1(0)- 1 )\cap \ker (C_2(0)- 1)\ \
\end{align*}
for a positive constant $c$.
\end{cor}

\bigskip

\begin{rem}
Corollary \ref{c:adn} and an elementary application of the
mini-max principle show that, in general, the operator $\Rr_R$ may
have exponentially decaying eigenvalues. Moreover, the number of
these eigenvalues is equal to
$$ \dim\big(\, \ker(C_1(0)-1)
\cap \ker(C_2(0)-1)\, \big) \ \ . $$ On the other hand, the
condition \eqref{Condition A}  and Remark \ref{r:c12} imply
\begin{equation}\label{e:inv12}
\ker(C_1(0)-1)\cap \ker(C_2(0)-1)=\{0\} \ \, ,
\end{equation}
hence it excludes the existence of exponentially small eigenvalues
of $\Rr_R$ under the condition \eqref{Condition A}. A simple
example where (\ref{e:inv12}) holds is the Dirac Laplacian over
the double of a manifold with boundary. It is easy to observe that
in this case we have $C_1(0)=-C_2(0)$ and there is no
exponentially small eigenvalues of $\Rr_R$ .
\end{rem}

\bigskip

Proposition \ref{p:prop-est2} suggests the introduction of the
operator $L(R)$ on $\ker(\Del_Y)$,
$$L(R) \ =\  \frac 1R
\biggm( \frac{\Id - C_1(0)}{2} + \frac{\Id - C_2(0)}{2}\biggm)  \
\ .$$

\bigskip

\begin{prop}\label{p:inv-1}
Assume that $\ker(C_1(0)-\Id)\cap \ker(C_2(0)-\Id)=\{0\}$ . Then
we have
\begin{equation}\label{e:inv-1}
\det  L(R) \ =\ R^{-h_Y} \det\big(\frac {\Id-C_{12}}2\big) \ \
\end{equation}
where $C_{12}:=C_1(0)\circ C_2(0)$.
\end{prop}

\begin{proof}
First of all, the assumption implies that the direct sum of the
ranges of the projections $\frac{\Id-C_1(0)}{2}$ ,
$\frac{\Id-C_2(0)}{2}$ spans the space $\ker(\Delta_Y)$ . It also
follows from the definition that we have a formula
$$
\det L(R)=R^{-h_Y} \det\biggm( \frac{\Id - C_1(0)}{2}+ \frac{\Id -
C_2(0)}{2} \biggm)  \ \ .
$$
Now, we use the fact that
\begin{equation}\label{e:inv212}
\frac {\Id - C_2(0)}2 = \biggm(\frac
{\Id-C_1(0)C_2(0)}2\biggm)^{-1} \frac {\Id+C_1(0)}2 \biggm(\frac
{\Id-C_1(0)C_2(0)}2\biggm) \ \, ,
\end{equation}
hence, essentially our concern is the determinant of the operator
acting on $\mathbb{C}^{h_Y}$ with the form
$$P + g^{-1}(\Id - P)g \ \ ,$$
putting $P=\frac{\Id - C_1(0)}{2} $ and $g=\frac
{\Id-C_1(0)C_2(0)}2$. We write
$$P + g^{-1}(\Id - P)g =
g^{-1}(gP + (\Id - P)g) \ \ .$$ The second operator on the right
side can be represented in the following form
\begin{equation}\label{e:3inv12}
gP + (\Id - P)g =
\begin{pmatrix}
PgP & 0 \\ 2(\Id - P)gP & (\Id - P)g(\Id - P)
\end{pmatrix}
\ \,
\end{equation} with respect to $\ran (P)\oplus
\ran(\mathrm{Id}-P)$. The corresponding decomposition for the
operator $P - g^{-1}(\Id - P)g$ is
$$g^{-1}\begin{pmatrix}
       PgP & 0 \\ 0 & -(\Id - P)g(\Id - P)
\end{pmatrix} \ \ .$$
This shows that
\begin{align*}
\det\biggm(& \frac{\Id-C_1(0)}{2}+ \frac{\Id-C_2(0)}{2} \biggm) =\
(-1)^{h_2}\det \biggm( \frac{\Id-C_1(0)}{2} -
\frac{\Id-C_{2}(0)}{2}\biggm)  \\
&=\ (-1)^{h_2} \det\biggm( \frac {\Id - C_{12}}2 \biggm)\det
C_2(0)
 = \ \det \biggm( \frac {\Id - C_{12}}2 \biggm) \ \,  .
\end{align*}
\end{proof}

\bigskip

{\bf Proof of Theorem \ref{t:DN0}}: Let $P^0$ and $P^{\perp}$
denote orthogonal projections onto the subspaces $\ker (\Del_Y)$
and $\ker( \Del_Y)^{\perp}$. For any trace class operator $L $
acting on $L^2(Y,E|_Y)$, we define
$$\Tr^0(L):=\Tr(P^0 L P^0) \quad ,\quad \Tr^{\perp}(L):=\Tr(P^{\perp} L
P^{\perp})\ \ .$$ We decompose $\Tr(e^{-t\Rr_R})$ into
$\Tr^0(e^{-t\Rr_R})$ and $\Tr^{\perp}(e^{-t\Rr_R})$ . By
Proposition \ref{p:prop-est2}, it is easy to see that the part
$\Tr^0(e^{-t\Rr_R})$ contributes by $\det L(R)$ up to the error of
the size $O(R^{-h_Y-1})$. By Proposition \ref{p:inv-1}, this is
$R^{-h_Y} \det(\frac {\Id-C_{12}}2)$ up to the error of the size
$O(R^{-h_Y-1})$.

\bigskip

Now let us see the contribution from $\Tr^{\perp}(e^{-t\Rr_R})$ .
Let us consider
\begin{align*}\label{e:res}
\ &\frac{i}{2\pi} \int_{\Gamma} \la^{-s}\ \Tr^{\perp}
((\Rr_R-\la)^{-1} - ( 2\sqrt{\Delta_Y} -\la)^{-1}) \ d\la\\
=\ & (-1)^k k! \frac{i}{2\pi} \int_{\Gamma} (s-1)^{-1}\cdots
(s-k)^{-k} \la^{-s+k}\\
&\qquad\qquad\qquad \quad \Tr^{\perp} ((\Rr_R-\la)^{-(k+1)} - (
2\sqrt{\Delta_Y} -\la)^{-(k+1)}) \ d\la \notag
\end{align*}
for sufficiently large $k$. Here, $\Gamma$ is a curve surrounding
$\{0\}\cup \mathbb{R}^-$ in $\mathbb{C}$. Let us remark that
$\Rr_R - 2\sqrt{\Delta_Y}$ is a smoothing operator. We refer the
proof of this fact to \cite{JPKW6}. Now the integrand on the right
side can be estimated as
\begin{align*}
& |\ \Tr^{\perp} ((\Rr_R-\la)^{-(k+1)} - ( 2\sqrt{\Delta_Y}
-\la)^{-(k+1)})\ |  \\
&\ \leq \   \frac{C}{|\la|^{k}+1} \ | \ \Tr^{\perp}( \Rr_R^{-1} -
(2\sqrt{\Delta_Y})^{-1})\ |
\end{align*}
for a positive constant $C$. Here $(2\sqrt{\Delta_Y})^{-1}$
denotes the inverse of $2\sqrt{\Delta_Y}$ over
$\ker(\Delta_Y)^{\perp}$. Now, we use Proposition \ref{p:sgrtd}
proved in Section 5, to show that the concerned integrand
converges to $0$ uniformly for every $s$ in the compact
neighborhood of $0$ as $R\to\infty$. Hence its derivative at $s=0$
converges to $0$ as $R\to\infty$. This completes the proof of
Theorem \ref{t:DN0} if we use
\begin{equation}\label{e:dmult}
\zd^*(2\sqrt{\Delta_Y}) = 2^{\z_{\Delta}(0)}\zd^* \sqrt{\Delta_Y}
\ \, .
\end{equation}

\bigskip

{\bf Proof of Corollary \ref{c:c9y}}: Let us now come back to the
BFK formula (\ref{e:BFK0}),
$$
\frac{\zd \Del_R}{\zd \Del_{1,R}\cdot \zd \Del_{2,R}}\ = \ C(Y)\,
\zd\Rr_R \ \ .
$$
We can use Theorem \ref{t:main thm} and Theorem \ref{t:DN0} to
find the exact value of the local constant $C(Y)$ . Let us recall
that $C(Y)$ does not depend on the adiabatic process. Now, we have
$$
2^{-h_Y} \sqrt{\zd^* \Del_Y} \cdot
\det\big(\frac{\mathrm{Id}-C_{12}}{2}\big) = \lim_{R\to\infty}
R^{h_Y} \frac{\zd \Del_R}{\zd \Del_{1,R}\cdot \zd \Del_{2,R}}
$$
$$= C(Y)\ \lim_{R\to\infty}
R^{h_Y}\zd\Rr_R = C(Y)\ 2^{\z_{\Delta_Y}(0)} \zd^*
\sqrt{\Delta_Y}\cdot \det\big(\frac{\mathrm{Id}-C_{12}}{2}\big) \
\ .$$ From this and the equality $\sqrt{\zd^* \Del_Y}=\zd^*
\sqrt{\Delta_Y}$, we conclude
\begin{equation*}
C(Y) = 2^{- \z_{\Delta_Y}(0) - h_Y} \ \, .
\end{equation*}

\vskip 1cm

\section{Proof of technical proposition}

\bigskip

In this section we present the proof of the following proposition,

\bigskip

\begin{prop}\label{p:sgrtd}
For $R\gg 0$, there exist positive constants $c_1$ and $c_2$ such
that
$$
| \ \Tr^{\perp}( \Rr_R^{-1} - (2\sqrt{\Delta_Y})^{-1} ) \ | \ \le
c_1e^{-c_2R^{\frac 12}} \ \ .
$$
\end{prop}

\bigskip

 Instead of using $2\sqrt{\Delta_Y}$, we
compare the operator $\Rr_R$ with the model operator $\Rr^c_R$ on
the cylinder defined as follows. We introduce the cylinder $N_R =
[-R,R] \times Y$ with the Laplacian $\Delta_R^c = -\partial_u^2 +
\Delta_Y$ subject to the Dirichlet boundary conditions at $\{\pm
R\}\times Y$. Now, we cut $N_R$ at $u=0$ and get the operator
$\Rr_R^c$ in an obvious way. An explicit computation shows that
the operator $\Rr^c_{R}$ converges to $2\sqrt{\Delta_Y}$
exponentially on the space $\ker (\Delta_Y)^{\perp}$ , more
precisely
$$|\ \Tr^{\perp}( \Rr^c_{R} -
2\sqrt{\Delta_Y})\ | \ \le \ c_3e^{-c_4R} \ \ $$ for some positive
constants $c_3,c_4$ . Therefore, it is sufficient to show
\begin{align}\label{e:me}
|\ \Tr^{\perp} ( \Rr_R^{-1} - (\Rr^c_R)^{-1} ) \ | \ \le
c_1e^{-c_2R^{\frac 12}} \ \ . \ \
\end{align}
In order to prove \eqref{e:me}, we recall the following formula
for $\Rr_R^{-1}$ established in \cite{BFK}, \cite{LP},
$$
\Rr_R^{-1} \ = \ \g \Delta_R^{-1} \g^*
$$
where $\g$ is the restriction map to $\{0\}\times Y$ and $\g^*$ is
the adjoint of $\g$. We combine this equality with
\begin{equation}\label{e:inv}
\Delta_R^{-1} \ = \ \int_0^{\infty} e^{-t\Delta_R}\ dt \ \ ,
\end{equation}
in order to reduce our problem to the heat kernel estimates. We
decompose the left side of (\ref{e:inv}) into two parts as follows
$$\int_0^{\infty} e^{-t\Delta_R}\ dt \ = \
\int_0^{R^{2-\e}} e^{-t\Delta_R}\ dt \ +\ \int_{R^{2-\e}}^{\infty}
e^{-t\Delta_R}\ dt \ \ .$$

\bigskip
We will consider the large and small time contributions separately
in the following lemmas.

\bigskip

\begin{lem} \label{l:est1} For $R\gg 0$,
there are positive constants $c_1$, $c_2$ such that
\[
|\, \Tr^{\perp}\big( \int_{R^{2-\e}}^{\infty}
 \g\, e^{-t\Delta_R} \g^* \ dt \big)\ |
\leq c_1e^{-c_2R^{1-\e}}
\]
and the same estimate holds for $\Delta^c_R$ .
\end{lem}

\begin{proof}
We note that
\begin{align}\label{e:lt}
\g\, e^{-t\Delta_R} \, \g^* =  \sum_k
e^{-t\la_k^2}\Phi_k(x)|_{u=0} \tensor \Phi_k^*(y)|_{u=0}
\end{align}
where $\{\la_k^2,\Phi_k\}$ is a spectral resolution of the
operator $\Delta_R$.  We split the restriction of the eigensection
$\Phi_k$ to $\{0\}\times Y$ into $\Phi_k^0$ the part in $\ker (
\Delta_Y)$ and ${\hat{\Phi}}_k$ the remaining part.  We employ an
argument similar to the proof of Lemma \ref{l:cut-est} to obtain
\begin{equation}\label{e:hcdn5}
\|{\hat{\Phi}}_k\| \ \le\ c_1e^{-\sqrt{\mu_{h_Y+1}^2 - \la_k^2}R}
\ \ .
\end{equation}
Here, we note that the right side of (\ref{e:hcdn5}) has to be
changed into the constant $c_1$ if $\la_k > \mu_{h_Y+1}$, and the
constant $c_1$ is independent of $k$ . We need to discuss only the
contribution determined by $\hat{\Phi}_k$ since we are concerning
only on $\Tr^{\perp}( \cdot)$. We split this contribution in
(\ref{e:lt}) into two parts, that is, the sums over all
eigenvalues $R^{-1} \leq \la_k^2 $  and  $ \la_k^2 < R^{-1}$.

 In order to discuss the
sum over the eigenvalues smaller than $R^{-1}$, we use
(\ref{e:hcdn5}) and the fact that each eigenvalue of $\Delta_R$ is
bounded from below by $\frac c{R^{2+ \frac {\e}2}}$ (since there
is no exponentially small eigenvalues). Then  we have
\begin{align}\label{e:lt3}
&\ \int_{R^{2-\e}}^{\infty} \biggm(\sum_{\la_k^2 < R^{-1}}
e^{-t\la_k^2} \| {\hat{\Phi}}_k \|^2\biggm) \, dt \ \leq  \ c_1
e^{- c_2 R} \int_{R^{2-\e}}^{\infty} \biggm(\sum_{\la_k^2 <
R^{-1}}e^{-t\la_k^2}\biggm) \, dt  \\
\leq & \ c_1 e^{- c_2 R}\, \Tr (e^{-{ \Delta}_R})\,
\int_{R^{2-\e}}^{\infty} e^{-{(t-1)}R^{-(2+\frac{\e}{2})}} \ dt
 \leq  \ c_{3} e^{- c_4 R} \notag
\end{align}
for positive constants $c_1, c_2, c_3, c_4$. We have used here the
obvious estimate
$$\Tr( e^{-\Delta_R}) \le c_5\, \text{vol}(M_R)
\le c_6R \ \ .$$ The sum over the eigenvalues $ R^{-1} \leq
\la^2_k$ can be estimated as
\begin{align}\label{e:lt1}
\int_{R^{2-\e}}^{\infty} \biggm( \sum_{R^{-1}\leq \la_k^2 }
e^{-t\la_k^2} \|{\hat{\Phi}}_k \|^2 \biggm)\, dt \  \leq \
c_1^2\int_{R^{2-\e}}^{\infty}  \biggm(\sum_{R^{-1}\leq \la_k^2}
e^{-t\la_k^2}\biggm)\, dt
\end{align}
$$\leq  \ c_1^2{\cdot}\Tr (
e^{-{\Delta}_R}) \int_{R^{2-\e}}^{\infty} e^{- \frac {t-1}R} \ dt\
\leq \ c_7R{\cdot}\int_{R^{2-\e}}^{\infty} e^{- \frac {t-1}R}\ dt
\ \leq \ c_8e^{-R^{1-\e}} \ \ .$$ The first claim follows from
(\ref{e:lt3}) and (\ref{e:lt1}). In the same way, we can show that
the same estimate holds for the operator $\Delta^c_{R}$.
\end{proof}

\bigskip

\begin{lem} \label{l:est2} For $R\gg0$,
there are positive constants $c_1,c_2$ such that
\begin{equation}\label{e:diff}
|\, \Tr^{\perp} \big(  \int^{R^{2-\e}}_0 \g\, (e^{-t\Delta_R} -
e^{-t\Delta^c_R})\, \g^* \ dt \big)\, | \leq c_1e^{-c_2R^{\e}} \ \
.
\end{equation}

\end{lem}

\begin{proof}
It is sufficient to show that the following term has the claimed
bound,
\[
\int_0^{R^{2-\e}} \int_Y \, \| \g\, (e^{-t\Delta_R}(x,x) -
e^{-t\Delta^c_R}(x,x))\, \g^* \| \, dy\ dt .
\]
For this, we apply \emph{Finite Propagation Speed Property for the
Wave Operator} to compare $\Delta_R$ over $M_R$ with $\Delta^c_R$
over $N_R$ where we identify the parts $N_{\frac R2}$ of these in
an obvious way. Then we obtain the estimate
\begin{equation*}
\| \Ee_R (t;x,y) - \Ee^c_R(t;x,y) \| \le c_3 e^{-c_4{\frac
{R^2}t}} \ \,
\end{equation*}
where $\Ee_R(t;x,y)$, $\Ee^c_R(t;x,y)$ are heat kernels of
$\Delta_R$, $\Delta^c_R$ respectively and $x, y \in N_{\frac R2}$.
Therefore, the following estimate holds
\begin{equation}\label{e:fps}
\| \g\, (e^{-t\Delta_R} - e^{-t\Delta^c_R})\, \g^* \| \le c_3
e^{-c_4{\frac {R^2}t}} \ \, .
\end{equation}
We combine (\ref{e:fps}) with the following inequality
\begin{align*}\label{e:las}
& c_3\int_0^{R^{2-\e}}e^{-c_4\frac {R^2}t}dt \  \leq \
c_1e^{-c_2R^{\e}} \ \ . \notag
\end{align*}
This completes the proof.

\end{proof}

\bigskip
Putting $\e=\frac 12$, Lemma \ref{l:est1} and \ref{l:est2}
complete the proof of Proposition \ref{p:sgrtd}.

\end{document}